\documentclass[a4paper]{article}
\usepackage{enumerate}
\usepackage[usenames, dvipsnames]{xcolor} %for colored text
\usepackage{amsfonts} %for \mathbb
\usepackage{dsfont} %for \mathds
\usepackage{amsmath} %for \text
\usepackage{amsthm} %for \proof
\usepackage{fullpage} %reduce margins
\usepackage{graphicx} % for \rotatebox

\usepackage{bm}

\newtheorem{theorem}{Theorem}[section]
\newtheorem{definition}{Definition}[section]
\newtheorem{lemma}{Lemma}[section]
\newtheorem{corollary}{Corollary}[section]
\newtheorem{proposition}{Proposition}[section]

\theoremstyle{definition}
\newtheorem{remark}{Remark}[section]
\newtheorem{example}{Example}[section]

\numberwithin{equation}{section}

\newcommand{\One}{\textbf{1}}

\newcommand{\E}{\ensuremath{\mathbb{E}}}
\newcommand{\R}{\ensuremath{\mathbb{R}}}

\newcommand{\Prob}{\ensuremath{\mathbb{P}}}

\newcommand{\N}{\ensuremath{\mathbb{N}}}

\newcommand{\Askript}{\ensuremath{\mathcal{A}}}
\newcommand{\Bskript}{\ensuremath{\mathcal{B}}}

\newcommand{\Tau}{\ensuremath{\mathcal{T}}}

\newcommand{\Ra}{\quad \Rightarrow \quad}

\newcommand{\BIGOP}[1]{\mathop{\mathchoice
{\raise-0.22em\hbox{\huge $#1$}}
{\raise-0.05em\hbox{\Large $#1$}}{\hbox{\large $#1$}}{#1}}}
\newcommand{\bigtimes}{\BIGOP{\times}}

\title{Embedded Markov chain approximations\\ in Skorokhod topologies}
\author{Bj\"orn B\"ottcher\footnote{\texttt{bjoern.boettcher@tu-dresden.de}, TU Dresden, Fakult\"at Mathematik, Institut f\"ur Math. Stochastik, 01062 Dresden, Germany}}
\date{}

\begin{document}
\maketitle
\begin{abstract}
We prove a $J_1$-tightness condition for embedded Markov chains and discuss the four Skorokhod topologies in a unified manner.

In order to approximate a continuous time stochastic process by discrete time Markov chains one has several options to embed the Markov chains into continuous time processes. On the one hand there is the Markov embedding, which uses exponential waiting times. On the other hand each Skorokhod topology naturally suggests a certain embedding. These are the step function embedding for $J_1$, the linear interpolation embedding for $M_1$, the multi step embedding for $J_2$ and a more general embedding for $M_2$. We show that the convergence of the step function embedding in $J_1$ implies the convergence of the other embeddings in the corresponding topologies. For the converse statement a $J_1$-tightness condition for embedded time-homogeneous Markov chains is given.

Additionally it is shown that $J_1$ convergence is equivalent to the joint convergence in $M_1$ and $J_2$. 
\end{abstract}

\noindent \textbf{Keywords}: {Markov chain embedding, tightness, Skorokhod space, Skorokhod topologies, jump processes, Markov chain approximation}                                     

\medskip
\noindent \textbf{MSC2010}: {Primary: 60B10 (Convergence of probability measures), 
60J75 (Jump processes), 
60J05 (Discrete-time Markov processes on general state spaces)
}   

\section{Introduction}

The space of right continuous functions with left limits plays a prominent role in the theory of stochastic processes. Skorokhod \cite{Skor1956} was the first to consider this space with various metrics. He introduced four topologies: $J_1, J_2, M_1$ and $M_2$. The main focus in the literature is on the $J_1$ topology (e.g. \cite{Bill1968,EthiKurt86,JacoShir2002}) and more recently on $M_1$ (e.g. \cite{AvraTaqq1989,Whit2002}). We will be concerned with all four. But note that there are further topologies on the Skorokhod space: e.g. the sequential topology of Jakubowski \cite{Jaku1997} and the pseudo-path topology by Meyer and Zheng \cite{MeyeZhen1984}.

Given the relations of Skorokhod's topologies, for a fixed sequence the convergence in a stronger topology implies the convergence in a weaker topology, i.e., $J_1$-convergence implies $M_1$ and $J_2$-convergence, and either of these implies $M_2$-convergence. But when one starts with discrete time processes there are many ways to embed these into continuous time processes, and most embeddings do not converge in all four topologies. Actually each of the four Skorokhod topologies suggests a particular embedding, the weaker the topology is the 'wilder' the embedding can be (see Section \ref{sec:embb}). Thus a natural question is: can we switch the topology and the corresponding embedding without losing convergence?

Consider Markov chains with time steps of size $\frac{1}{n}$ and let $n$ tend to infinity. In order to discuss a continuous time limit it is necessary to embed the chains into continuous time processes. In our general setting the limit can be a process with jumps. For processes with continuous paths Sato \cite{Sato1977} discussed a closely related problem: he showed that linearly interpolated Markov chains converge with respect to the uniform topology (in the space of continuous functions) if and only if the step function embedded Markov chains converge to a continuous process with respect to the $J_1$ topology (in the Skorokhod space). Our result allows in particular, cf. Example \ref{rem-feller}, to extend Markov chain approximations for Feller processes (cf. \cite{BoetSchi2009,BoetSchiWang2013}) to different embeddings. More general, we provide a $J_1$-tightness condition for Markov chains, see Theorem \ref{J1-mc-convergence}. 

It turns out that, in the above setting, convergence is always preserved when switching from a topology to a weaker topology (and to the corresponding embedding), see Corollary \ref{simpleimplication} and Theorem \ref{embedding-equiv}. For the converse direction naturally some additional assumption is needed, see (Counter-)Examples \ref{counter-ex} and Corollary \ref{conv-tight-equi}. 

In the next section we introduce the Skorohod space and the topologies $J_1, J_2,$  $M_1$ and $M_2$ in a unified framework, which consolidates the literature e.g. \cite{Skor1956,Whit2002,Poma1976}. In particular we recall their relations and several representations. The relation between $J_1$ and the combination of $J_2$ and $M_1$ (Lemma \ref{osc-relations}) seems to be neglected in the literature. It goes back to a remark without proof of Skorokhod \cite[2.2.10-13]{Skor1956}. In Section \ref{sec:embb} the embeddings are introduced and their relations are discussed. In Section \ref{sec:mc} a $J_1$-tightness condition (Theorem \ref{J1-mc-convergence}) for embedded Markov chains is presented, it enables us to switch from a weaker to a stronger topology (and to the corresponding embedding; see Corollary \ref{conv-tight-equi}). The paper closes with the proof of Lemma  \ref{osc-relations}.

\section{The Skorokhod space and its topologies}
Throughout the paper segments between points $x,y \in \R^d$ are denoted by 
$$[[x,y]]:=\{z\in\R^d\,|\,z = \alpha x + (1-\alpha) y \text{ for some } \alpha \in [0,1] \},$$
and $\|.\|_\infty$ denotes the supremum norm. Limits without superscript, e.g. $f_n \to f$, are meant in the Euclidean distance. Unless stated otherwise, limits are considered for the index tending to infinity, e.g. $n\to \infty$, and the dimension $d\in\N$ is arbitrary.

\begin{definition}
A function $f:[0,1] \to \R^d$ is \textbf{right continuous with left limits} (rcll)  if 
$$f(s+):=\lim_{t\downarrow s} f(t) = f(s) \text{ and } f(s-):=\lim_{t\uparrow s} f(t) \text{ exist in $\R^d$.}$$
The \textbf{Skorokhod space} is
$$D[0,1]:=D_d[0,1]:= \left\{f:[0,1]\to \R^d \,\big|\, f \text{ is rcll, and $f$ is left continuous at $1$}\right\}.$$
\end{definition}

On the Skorokhod space several metrics can be defined.
\begin{definition}
Let $f,f_1,f_2 \in D[0,1].$\\
The \textbf{incomplete graph} of $f$ is $$\Gamma_f:= \{(z,t)\in \R^d\times[0,1] \,|\, z= f(t-) \text{ or } z= f(t)\}.$$
The \textbf{complete graph} of $f$ is 
$$\overline{\Gamma_f}:= \{(z,t)\in \R^d\times[0,1] \,|\, z \in [[f(t-),f(t)]]\}.$$
An order is defined on $\overline{\Gamma_f}$ by
$$(z_1,t_1) \leq (z_2,t_2) \text{ if }\begin{array}{l}\text{either $t_1< t_2,$} \\\text{or $t_1=t_2, |f(t_1-)-z_1|\leq |f(t_2-)-z_2|$} \end{array}$$
and the families of \textbf{parametric representations} of $\overline{\Gamma_f}$ are given by
\begin{equation*}
\begin{split}
\pi(\overline{\Gamma_f}) = \{(u,r)\, |&\, u:[0,1]\to \R^d \text{ continuous}, r:[0,1]\to [0,1]\text{ continuous}, \\
&(u,r)\text{ is non decreasing and } (u,r)[0,1] = \overline{\Gamma_f}\},
\end{split}
\end{equation*}
\begin{equation*}
\begin{split}
{\tilde \pi(\overline{\Gamma_f})} = \{(u,r)\, |\, &u:[0,1]\to \R^d \text{ continuous}, r:[0,1]\to [0,1] \text{ continuous},\\
& {r\text{ is non decreasing and }} (u,r)[0,1] = \overline{\Gamma_f}\}.
\end{split}
\end{equation*}
The sets of \textbf{time transformations} are $\Lambda := \{\lambda:[0,1]\to [0,1]\,|\,\lambda \text{ is bijective}\}$ and $\Lambda_c:= \Lambda\cap C[0,1].$ Let $\mathrm{id}$ denote the identity function on $[0,1]$ and then \textbf{metrics on Skorokhod space} are given by\\
\begin{equation*}
\begin{aligned}
\bm{ d_{J_1}}(f_1,f_2) &:= \inf_{\lambda\in \Lambda_c} \{\|f_1\circ \lambda - f_2\|_\infty \lor \|\lambda - \mathrm{id}\|_\infty\},\\
{\tilde d_{J_2}}(f_1,f_2) &:= \inf_{{\lambda\in \Lambda_{\phantom{c}}}} \{\|f_1\circ \lambda - f_2\|_\infty \lor \|\lambda - \mathrm{id}\|_\infty\},\\
\bm{ d_{M_1}}(f_1,f_2) &:= \inf_{\substack{(u_j,r_j)\in\pi(\overline{\Gamma_{f_j}})\\ j=1,2}}\{\|u_1-u_2\|_\infty \lor \|r_1-r_2\|_\infty\},\\
{\tilde d_{M_2}}(f_1,f_2) &:= \inf_{\substack{(u_j,r_j)\in{\tilde \pi(\overline{\Gamma_{f_j}})}\\ j=1,2}}\{\|u_1-u_2\|_\infty \lor \|r_1-r_2\|_\infty\},\\
\bm{ d_{J_2}}(f_1,f_2) &:= m_H(\Gamma_{f_1},\Gamma_{f_2}),\\
\bm{ d_{M_2}}(f_1,f_2) &:= m_H(\overline{\Gamma_{f_1}},\overline{\Gamma_{f_2}})
\end{aligned}
\end{equation*}
where $m_H(A,B) = \inf\{\varepsilon>0 : A\subset B^\varepsilon, B\subset A^\varepsilon\}$ is the Hausdorff distance with 

$B^\varepsilon:=\{(x,t)\in\R^d\times[0,1]\, |\, \exists (y,s) \in B : |x-y|\lor|t-s|\leq \varepsilon\}.$

\end{definition}
In the following $\Tau$ will always denote one of $J_1,J_2,M_1,M_2$. Note that $d_\Tau$ is a metric on $D[0,1]$. This is clear by the definition for $J_2$ and $M_2$, for $J_1$ see \cite{Bill1968} and for $M_1$ see \cite{Whit2002}. 
For $f,f_n\in D[0,1]$ $(n\in\N)$ we say \textbf{$f_n$ converges to $f$ in $\Tau$} (in $\Tau$-topology), if $d_\Tau(f_n,f) \to 0.$ The convergence is denoted by
$$f_n \xrightarrow{\Tau} f.$$

\begin{remark} \label{j2m2-orig}
Skorokhod \cite{Skor1956} used $\tilde d_{J_2}, \tilde d_{M_2}$ to introduce the topologies $J_2$ and $M_2$, respectively. But $d_{J_2}, d_{M_2}$ yield the same topologies \cite[II.4.1 p.\ 82, II.4.2 p.\ 83]{Poma1976}. Note that on $\R^d$ with $d>1$ one could also define a complete graph by including for each coordinate the whole interval between the start and endpoints of the jumps, i.e. $\bigtimes_{i=1}^d [[f_i(t-),f_i(t)]]$ where $f_i$ denotes the $i$th component of $f$. The definition above only uses the linear interpolation, i.e., $[[f(t-),f(t)]]$, thus here the $M$ topologies are strong in the sense of Whitt (\cite{Whit2002} Section 12.3 and 12.10). 

Regarding completeness of the corresponding metric spaces see Remark \ref{rem-fdd-d}.1.
\end{remark}

The convergence in these topologies can also be characterized by oscillation functions.

\begin{definition}(Oscillation functions) 

Define for $x,x_1,x_2\in\R^d$
\begin{equation}
\begin{aligned}
J(x,x_1,x_2) &:= |x-x_1|\land |x-x_2|,\\
M(x,x_1,x_2) &:= \big|x-[[x_1,x_2]]\big| := \inf_{y\in [[x_1,x_2]]}|x-y|
\end{aligned}
\end{equation}
and for $\delta >0$
$$T_1(\delta) := \{(t,t_1,t_2)\ |\ (t-\delta)\lor 0 \leq t_1 < t< t_2\leq (t+\delta)\land 1\},$$
\begin{equation*}
\begin{split}
T_2(\delta) := \{ (t,t_1,t_2)\ |\ t\in[0,1], &t_1\in[(t-\delta)\lor 0, (t-\delta)\lor 0 + \frac{\delta}{2}],\\ 
&t_2\in[(t+\delta)\land 1 - \frac{\delta}{2}, (t+\delta)\land 1]\}.
\end{split}
\end{equation*}
The \textbf{oscillation functions} for $f:[0,1]\to \R^d$, $\delta >0$ and $i=1,2$ are
\begin{equation}
\begin{aligned}
\Delta_{J_i}(\delta,f) &:= \sup_{(t,t_1,t_2)\in T_i(\delta)} J(f(t),f(t_1),f(t_2)),\\
\Delta_{M_i}(\delta,f) &:= \sup_{(t,t_1,t_2)\in T_i(\delta)} M(f(t),f(t_1),f(t_2)),\\
\Delta_U^{\{0,1\}}(\delta,f) &:= \sup_{0<t< \delta } |f(0)- f(t)| + \sup_{1-\delta< t <1} |f(1)-f(t)|.
\end{aligned}
\end{equation}
\end{definition}

The following theorem states the fundamental relation of the oscillation functions and the metrics.

\begin{theorem}\label{thm-oscill-char} Let $f_n, f \in D[0,1].$ Then
$$
f_n \xrightarrow{\Tau} f \ \Leftrightarrow \begin{array}{l} \text{i)\phantom{i}} f_n(t)\to f(t), \forall t \in B, \text{ where }0,1 \in B \text{ and}\\ 
\phantom{\text{ii)} f_n(t)\to f(t), \forall t \in B, } \text{ $B$ is a dense subset of }[0,1],\\
\text{ii)}\lim_{\delta \downarrow 0} \limsup_{n\to\infty} \Delta_\Tau(\delta, f_n) = 0. \end{array}
$$
\end{theorem}
\begin{proof}
The proofs can be found for $M_2$ in  \cite[2.3.4]{Skor1956}, for $M_1$ in \cite[2.4.1]{Skor1956}, for $J_2$ in \cite[II.4.4]{Poma1976},  \cite[2.5.3]{Skor1956} and for $J_1$ in \cite[2.6.1]{Skor1956}.
\end{proof}

Note that the oscillation functions satisfy the following relations.

\begin{lemma} \label{osc-relations} Let $\delta >0$ and $f:[0,1] \to \R^d.$ Then
\begin{equation}
\Delta_{M_2}(\delta,f)\ \begin{array}{r} \mathbin{\raisebox{.2em}{\rotatebox[origin=c]{20}{$\leq$}}} \\ \mathbin{\rotatebox[origin=c]{-20}{$\leq$}}  \end{array} \begin{array}{c}\Delta_{J_2}(\delta,f)\\\\\Delta_{M_1}(\delta,f) \end{array} \begin{array}{l} \mathbin{\raisebox{.2em}{\rotatebox[origin=c]{-20}{$\leq$}}} \\ \mathbin{\rotatebox[origin=c]{20}{$\leq$}}  \end{array} \ \Delta_{J_1}(\delta,f)\  \leq \ \Delta_{M_1}(\delta,f) + \Delta_{J_2}(\delta,f).
\end{equation}
\end{lemma}
\begin{proof}
The first four inequalities follow directly from the definition of the oscillation functions, since $M(x,x_1,x_2) \leq J(x,x_1,x_2)$ and $T_2(\delta) \subset T_1(\delta).$ The last inequality is proved in Section \ref{sec:proofs}.
\end{proof}

Thus we have the following relations of the convergences

\begin{equation}
\label{convergence-relations}
M_2\ \begin{array}{r} \mathbin{\raisebox{.2em}{\rotatebox[origin=c]{20}{$\Leftarrow$}}} \\ \mathbin{\rotatebox[origin=c]{-20}{$\Leftarrow$}}  \end{array} \begin{array}{c}J_2\\\\M_1 \end{array} \begin{array}{l} \mathbin{\raisebox{.2em}{\rotatebox[origin=c]{-20}{$\Leftarrow$}}} \\ \mathbin{\rotatebox[origin=c]{20}{$\Leftarrow$}}  \end{array} \ J_1\  \Leftrightarrow \ M_1 + J_2.
\end{equation}

As remarked by Skorokhod \cite[2.2.10-13]{Skor1956} there are further equivalent characterizations of the convergence in these topologies for functions in $\R^d$ with $d=1$.
\begin{theorem}\label{convergence-characterization} Let $f_n, f\in D_1[0,1]$.  
\begin{enumerate} 
\item $M_2$ is characterized by the convergence of the local extrema: 
$$f_n \xrightarrow{M_2} f \quad\Leftrightarrow\quad \inf_{t\in[t_1,t_2]}f_n(t) \to \inf_{t\in[t_1,t_2]}f(t)\text{ and } \sup_{t\in[t_1,t_2]} f_n(t) \to \sup_{t\in[t_1,t_2]} f(f) $$
for all $t_1,t_2$ being points of continuity of $f$.
\item $M_1$ is characterized by the convergence of the number of oscillations:
$$f_n \xrightarrow{M_1} f \quad\Leftrightarrow \quad\nu_{[t_1,t_2]}^{[a,b]}f_n \to \nu_{[t_1,t_2]}^{[a,b]}f$$
for all $t_1,t_2$ being points of continuity of $f$ and almost all $a<b$. Here $\nu_{[t_1,t_2]}^{[a,b]}f$ is the largest $k$ such that there exist $t^{(0)} < \ldots < t^{(k)}$ in $[t_1,t_2]$ with $f(t^{(0)})\leq a,$ $f(t^{(1)})\geq b,$ $f(t^{(2)})\leq a,\ldots$
\item $J_2$ is characterized by the convergence of the first overshoots:
$$f_n \xrightarrow{J_2} f \quad\Leftrightarrow \quad\gamma^+_{[t_1,t_2],a}f_n \to \gamma^+_{[t_1,t_2],a}f$$
for all $t_1,t_2$ being points of continuity of $f$ and almost all $a$. Here, using the convention $\inf \emptyset = 1$, define $\tau_{a,f} := \inf\{t\in[0,1]\,|\, f(t)\geq a\}$ and
$$\gamma^+_{[0,1],a}f := \begin{cases} f(\tau_{a,f}) -a &, \tau_{a,f}< 1,\\ -1&, \text{ otherwise,}\end{cases}$$
and in general use $\gamma_{[t_1,t_2],a}^+f:=\gamma_{[0,1],a}^+\tilde f$ with $\tilde f(t):=\begin{cases} f(t_1+) &, t\leq t_1,\\ f(t) &, t\in(t_1,t_2),\\ f(t_2-) &, t\geq t_2.\end{cases}$ 
\item $J_1$ is characterized by the convergence of the first overshoots and the number of oscillations:
$$f_n \xrightarrow{J_1} f \quad\Leftrightarrow\quad  \gamma^+_{[t_1,t_2],a}f_n \to \gamma^+_{[t_1,t_2],a}f \text{ and } \nu_{[t_1,t_2]}^{[a,b]}f_n \to \nu_{[t_1,t_2]}^{[a,b]}f $$
for all $t_1,t_2$ being points of continuity of $f$ and almost all $a<b$. For the definition of $\gamma$ and $\nu$ see 2. and 3.
\end{enumerate}
\end{theorem}
\begin{proof}
The first and third statement are a consequence of the definition of these metrics via the Hausdorff metric. The second statement can be found in Whitt \cite[Thm. 12.7.4, p.\ 412]{Whit2002}. The last statement is due to the equivalence of the convergences (Lemma \ref{osc-relations}): $J_1 \Leftrightarrow M_1+J_2.$
\end{proof}
\begin{remark}\label{extension-of-01}
\begin{enumerate}
\item The characterizations in Theorem \ref{convergence-characterization} are tailored to $d=1$. For higher dimensions Whitt \cite[Theorem 12.7.2]{Whit2002} showed, for example, that 
\begin{equation}\label{m1-dim-d-1}
f_n\xrightarrow{M_1}f  \quad \Leftrightarrow \quad \eta\cdot f_n \xrightarrow{M_1} \eta\cdot f \text{ as functions in $D_1[0,1]$ for all }\eta \in \R^d
\end{equation}
for $f_n, f\in D_d[0,1].$
\item Throughout this section we only considered $D[0,1]$. By replacing $1$ by $T\in (0,\infty)$ we have an obvious extension to $D[0,T]$. An approach to define convergence for  $f_n,f\in D[0,\infty)$ is
\begin{equation}
\begin{split}
&f_n\xrightarrow{\Tau}f  \\
\Leftrightarrow \quad &f_n\big|_{[0,T]} \xrightarrow{\Tau} f\big|_{[0,T]} \text{ in $D[0,T]$ for all } T \in \{t\,|\,f  \text{ is continuous in } t\}.
\end{split}
\end{equation}
For further details on the extension to $D[0,\infty)$ see Lindvall \cite{Lind1973}.
\end{enumerate}
\end{remark}

\section{Embeddings and approximations}\label{sec:embb}
Let $n\in\N$ and $y^{(n)}$ be a sequence $(y_0^{(n)}, y_1^{(n)}, y_2^{(n)},\ldots)$ in $\R^d.$ Define for each topology $\Tau$ the embeddings $x^{n,\Tau}$ to be functions in $D[0,1]$ such that for $k<n$
\begin{equation}
x^{n,\Tau}(\textstyle\frac{k}{n}) = y_k^{(n)}
\end{equation}
and for all $t\in (\frac{k}{n},  \frac{k+1}{n})$
\begin{equation} \label{def-embedding}
\begin{aligned}
x^{n,J_1}(t)&=y^{(n)}_k, &\text{ -- step functions for $J_1$ --}\\
x^{n,M_1}(t)&= y_{k}^{(n)} + (t-\textstyle\frac{k}{n}) (y_{k+1}^{(n)}-y_k^{(n)}), &\text{ -- linear interpolation for $M_1$ --}\\
x^{n,J_2}(t)&\in\{y_k^{(n)},y_{k+1}^{(n)}\}, &\text{ -- multiple steps for $J_2$ --}\\
x^{n,M_2}(t)&\in[[y_k^{(n)},y_{k+1}^{(n)}]]. &\text{ -- any rcll function for $M_2$ --}
\end{aligned}
\end{equation}
Note that the requirement $x^{n,\Tau} \in D[0,1]$ ensures that $x^{n,\Tau}(1) = \lim_{t\nearrow 1} x^{n,\Tau}(t).$ Clearly  in the above definition only those $k$ with $k<n$ are used, but in the next section it will be convenient that each $y^{(n)}$ is a countable sequence.

Lemma \ref{osc-relations} implies the following result.
\begin{corollary} \label{simpleimplication}
$x^{n,J_1}$ converges in $J_1$ implies that $x^{n,J_1}$ also converges in $\Tau$.
\end{corollary}

Moreover in a given topology we can always switch between its embedding and the $J_1$ embedding.
\begin{theorem} \label{embedding-equiv}
$x^{n,\Tau}$ converges in $\Tau$ if and only if $x^{n,J_1}$ converges in $\Tau$.
\end{theorem}
\begin{proof}
By the definition of the metrics and the embeddings
\begin{equation}
d_{J_2}(x^{n,J_2},x^{n,J_1}) \leq \frac{1}{n} \quad \text{ and } \quad d_{M_2}(x^{n,M_2},x^{n,J_1}) \leq \frac{1}{n}.
\end{equation}
For $M_1$ we use \eqref{m1-dim-d-1} and Theorem \ref{convergence-characterization}: Note that for all $\eta\in\R^d$ and all $t_1,t_2$ which are points of continuity of the limit and almost all $a<b$ 
\begin{equation}
\left|\nu_{[t_1,t_2]}^{[a,b]} (\eta \cdot x^{n,M_1}) - \nu_{[t_1,t_2]}^{[a,b]}(\eta \cdot x^{n,J_1})\right| \xrightarrow{n\to \infty} 0,
\end{equation}
since for $\frac{k-1}{n} < t_1 \leq \frac{k}{n}$ and $\frac{l}{n} \leq t_2 < \frac{l+1}{n}$ the number of oscillations coincides for the segment from $\frac{k}{n}$ to $\frac{l}{n}$. In the limit no overshoot appears at the two boundary segments since $t_1$ and $t_2$ are points of continuity.  
Thus if the limit is in $D[0,1]$ the statement follows by the triangle inequality.
\end{proof}

We close this section with basic counterexamples which show that the converse implication of Corollary \ref{simpleimplication} fails.

\begin{example} \label{counter-ex}
Let $n\geq 4$ and $k\in \N\cup\{0\}.$

\begin{enumerate}
\item Let $y_k^{(n)} = 0$ for $\frac{k}{n}< \frac{1}{2}$, $y_k^{(n)}=\frac{1}{2}$ for $\frac{k-1}{n}< \frac{1}{2} \leq \frac{k}{n}$ and $y_k^{(n)}=1$ otherwise. Then $x^{n,J_1}$ converges to $\One_{[\frac{1}{2},1]}$ in $M_1$, but not in $J_1$ and not in $J_2$.
\item Let $y_k^{(n)} = 0$ for $\frac{k}{n}< \frac{1}{2}$, $y_k^{(n)}=1$ for $\frac{k-1}{n}< \frac{1}{2} \leq \frac{k}{n}$, $y_k^{(n)}=0$ for $\frac{k-2}{n}< \frac{1}{2} \leq \frac{k-1}{n}$ and $y_k^{(n)}=1$ otherwise. Then $x^{n,J_1}$ converges to $\One_{[\frac{1}{2},1]}$ in $J_2$, but not in $J_1$ and not in $M_1$.
\item Let $y_k^{(n)} = 0$ for $\frac{k}{n}< \frac{1}{2}$, $y_k^{(n)}=\frac{1}{2}$ for $\frac{k-1}{n}< \frac{1}{2} \leq \frac{k}{n}$, $y_k^{(n)}=1$ for $\frac{k-2}{n}< \frac{1}{2} \leq \frac{k-1}{n}$,  $y_k^{(n)}=0$ for $\frac{k-3}{n}< \frac{1}{2} \leq \frac{k-2}{n}$ and $y_k^{(n)}=1$ otherwise. Then $x^{n,J_1}$ converges to $\One_{[\frac{1}{2},1]}$ in $M_2$, but not in $J_2$ and not in $M_1$, and thus not in $J_1$.
\end{enumerate} 
\end{example}

\section{Convergence of processes and Markov chains} \label{sec:mc}

Let $X,X^{(n)}$ $(n\in\N)$ be $D[0,1]$-valued random variables on some probability space $(\Omega,\Askript,\Prob).$ To fix notations we recall the following standard definitions.

\begin{definition} 
\begin{enumerate}
\item $X^{(n)} \xrightarrow{d} X \text{ w.r.t. $\Tau:$ } \E(G(X^{(n)}))\to \E(G(X))$ for all bounded and $\Tau$-continuous functions $G:D[0,1]\to \R,$
\item $X_t^{(n)} \xrightarrow{d} X_t:$ $\E(g(X_t^{(n)}))\to \E(g(X_t))$ for all bounded and continuous functions $g:\R^d\to \R^d,$ 
\item $X^{(n)} \xrightarrow{\Prob} X \text{ w.r.t. } \Tau:$  $\lim_{n\to\infty}\Prob(d_\Tau(X^{(n)},X)>\varepsilon) = 0$ for all $\varepsilon>0,$ 
\item $X^{(n)} \xrightarrow{\text{fdd}} X$ on $I:$ $(X_{t_1}^{(n)},\ldots, X_{t_k}^{(n)}) \xrightarrow{d} (X_{t_1},\ldots, X_{t_k})$ for all $t_i\in I$, 
\item $(X^{(n)})_{n\in\N}$ is $\Tau$-tight: for all $\varepsilon>0$ exists a $\Tau$-compact set $K\in\Bskript(D[0,1])$ such that $\sup_n\Prob(X^{(n)}\in K^c) \leq \varepsilon$,
\item $(X^{(n)})_{n\in\N}$ is relative $\Tau$-compact: for every subsequence $(X^{(n_{k})})_{k\in\N}$ exists a further subsequence $(X^{(n_{k_l})})_{l\in\N}$ and a $D[0,1]$-valued random variable $Y$ such that $X^{(n_{k_l})} \xrightarrow{d} Y$ w.r.t. $\Tau$. 
\end{enumerate}
\end{definition}

The following result is the standard tool to handle convergence on $D[0,1]$. We include a sketch of the proof since we are going to point out a particular detail later.

\begin{theorem} \label{thm-fdd-d}
$$ X^{(n)}\xrightarrow{d} X\text{ w.r.t. }\Tau \quad \Leftrightarrow \quad \begin{array}{l} \text{i)\phantom{i} }X^{(n)} \xrightarrow{\text{fdd}} X\text{ on a dense subset of }[0,1],\\ \text{ii) } (X^{(n)})_{n\in\N}\text{ is relatively $\Tau$-compact.}\end{array}$$
\end{theorem}
\begin{proof}
'$\Leftarrow$': By ii) every subsequence of $X^{(n)}$ has a converging subsequence whose limit has by i) the same finite dimensional distributions as $X$. The finite dimensional distributions define uniquely the distribution of a process in $D[0,1]$, thus the limit is $X$.

'$\Rightarrow$': The projection $\pi_t:D[0,1] \to \R^d,$ $\pi_t(f):= f(t)$ is measurable, $X^{(n)} \xrightarrow{d} X$ and the set $T:=\{t: \Prob(|X_{t}-X_{t-}|)>0\}$ is countable. Thus for all $t\in T^c$ 
$$\Prob(X\in \{f\in D[0,1] : \pi_t(f) \text{ is discontinous at t}\}) = \Prob(|X_{t}-X_{t-}|>0)=0$$
and the statement follows by a continuous mapping theorem, e.g. \cite[Theorem 3.4.3]{Whit2002}. 
\end{proof}

\begin{remark}\label{rem-fdd-d}
\begin{enumerate}
\item A sufficient condition for relative $\Tau$-compactness is given by Prohorov's Theorem:
$$ \Tau\text{-tightness} \Ra \text{relative }\Tau\text{-compactness}.$$
The converse holds if $D[0,1]$ is, with the topology induced by $d_\Tau$, a complete and separable space. This is the case for $J_1$ and $M_1$, but for $J_2$ and $M_2$ it is still an open problem. To avoid confusion, note that the metric spaces $(D[0,1],d_{J_1})$ and $(D[0,1],d_{M_1})$ are not complete - but there exist complete metrics which generate the same topologies (one approach to construct these complete metrics is to add to the given metrics the L\'evy distance of distributions obtained via the oscillation functions; see Section 12.8 in \cite{Whit2002}).
\item Note that we assumed that $X$ is $D[0,1]$-valued. For condtion i) in Theorem \ref{thm-fdd-d} this can be relaxed, at least if $J_1$ is considered. Topsoe \cite[Theorem 2]{Tops1969} showed that $J_1$-tightness and the mere convergence of the finite dimensional distributions on a dense subset are sufficient to identify a process in $D[0,1]$ which is the $J_1$ limit.
\item Looking at the proof of Theorem \ref{thm-fdd-d} (see also \cite[Theorem 3.14]{JacoShir2002}) note that if the process $X$ is stochastically continuous, i.e.,
\begin{equation}
\label{equ-stoch-cont}
\forall s\in[0,1]\ \forall \varepsilon>0:\ \lim_{t\to s}\Prob(|X_t-X_s|>\varepsilon)=0,
\end{equation}
then the dense subset of $[0,1]$ can be taken to be the whole set $[0,1].$ See also Proposition \ref{stoch-continuity} below. 
\item A necessary condition for $ X^n\xrightarrow{d} X\text{ w.r.t. }\Tau$ is
\begin{equation}\label{global-bound}
\forall \varepsilon>0\ \exists R>0:\ \sup_{n}\Prob(\|X^{(n)}\|_\infty \geq R) <\varepsilon,
\end{equation}
since otherwise some mass would dissipate and hence $X$ would have, with positive probability, values not in $D[0,1].$ Also note that $\{f\in D[0,1]\,|\, \|f\|_\infty <R\}$ is not $\Tau$-compact. Since e.g. $\left(\One_{[\frac{1}{2},\frac{1}{2}+\frac{1}{n})}\right)_{n\geq 2}$ has no converging subsequence.
\end{enumerate}
\end{remark}

For $J_1$ there are several conditions for tightness, we will start with a standard result (see e.g. \cite[Theorem 3.21]{JacoShir2002}).

\begin{theorem}
\label{J1-tight}
The sequence of processes $(X^{(n)})_{n\in\N}$ is $J_1$-tight if and only if \eqref{global-bound} holds and 
\begin{equation}
\label{tight-oszill}
\forall \varepsilon>0:\ \lim_{\delta\downarrow 0} \sup_n \Prob\left(\Delta_{J_1}(\delta, X^{(n)})+\Delta_U^{\{0,1\}}(\delta, X^{(n)})>\varepsilon\right) = 0.
\end{equation}
\end{theorem}
\begin{remark}
For $M_1$ a result analogous to Theorem \ref{J1-tight} holds (cf. \cite[Theorem 12.12.3]{Whit2002}). For $M_2$ and $J_2$ the corresponding version of condition \eqref{tight-oszill} is not necessary for compactness (cf. \cite[2.7.2-4]{Skor1956}).
\end{remark}

A well known sufficient $J_1$-tightness condition is due to Aldous \cite{Aldo1978}. 
\begin{theorem}[Aldous \cite{Aldo1978}] \label{thm-aldo}
The sequence $(X^{(n)})_{n\in\N}$ is $J_1$-tight if 
\begin{equation} \label{aldous-cond}
\forall \varepsilon>0:\ \lim_{n\to \infty}\Prob\left(|X^{(n)}_{(\tau_n+t_n)\land 1}-X_{\tau_n}^{(n)}|>\varepsilon\right) = 0
\end{equation}
for all sequences $(\tau_n)_{n\in\N}$, with $\tau_n$ being a stopping time for $X^{(n)}$, and all sequences $(t_n)_{n\in\N}$ with $t_n\geq 0,\, t_n\to 0.$
\end{theorem}

As a motivation we also recall a closely related result by Gikhman and Skorokhod \cite[Theorem 4, p. 431]{GikhSkor1974} formulated here for the time homogeneous setting:

\begin{theorem}[Gikhman, Skorokhod \cite{GikhSkor1974}] \label{thm-gisko}
Let $X, X^{(n)}$ be time homogeneous Markov processes with $X^{(n)} \xrightarrow{\text{fdd}} X$ on some dense subset of $[0,1]$ and 
\begin{equation} \label{cond-gisko}
\forall\varepsilon>0:\ \lim_{h\downarrow 0} \limsup_{n\to \infty} \sup_{\substack{x\in\R^d\\t\leq h}} \Prob\left(|X_t^{(n)}-x|>\varepsilon\, \bigg|\, X_0^{(n)}=x\right) =0,
\end{equation}
then $X^{(n)} \xrightarrow{d} X$ w.r.t. $J_1.$
\end{theorem}

Thus \eqref{cond-gisko} is a $J_1$-tightness condition, actually ensuring that the limit is spatial-uniformly stochastically continuous from the right. Aldous tightness condition \eqref{aldous-cond} and condition \eqref{cond-gisko} are both not necessary for convergence, a counterexample is a process with a fixed jump, e.g. consider the deterministic time homogeneous Markov process whose transition probabilities for $t>0$ and $x\in \R$ are
\begin{equation*}
\begin{aligned}
\Prob(X_t = x \,|\, X_0 = x) = 1& &\text{ for all } x\in [0,1), \\
\Prob(X_t = x+t \,|\, X_0 = x) = 1& &\text{ for all } x\in [1,\infty)\text{ or } x\in (-\infty,0\land (-t)),\\ 
\Prob(X_t = x+t+1 \,|\, X_0 = x) = 1& &\text{ for all } x\in (-\infty,0)\text{ and } t+x \geq 0. \\
\end{aligned}
\end{equation*}

Incidentally, this counterexample also shows that for time homogeneous Markov processes stochastic continuity \eqref{equ-stoch-cont} is stronger than stochastic continuity from the right, i.e., $\lim_{t\downarrow 0}\Prob(|X_t-x|>\varepsilon \ |\ X_0=x) = 0$ for all $\varepsilon>0$ and $x\in\R^d.$ 
In fact the following holds.
\begin{proposition}\label{stoch-continuity}
Let $X$ be a $D[0,1]$ valued process then the following are equivalent:
\begin{enumerate} 
\item $X$ is stochastically continuous: 
$$\forall s\in[0,1]\ \forall \varepsilon>0:\ \lim_{t\to s}\Prob(|X_t-X_s|>\varepsilon)=0,$$
\item $X$ has no fixed discontinuities: $\forall s\in[0,1]: \ \Prob(|X_s-X_{s-}|>0)=0.$
\end{enumerate}
If additionally $X$ is a time homogeneous Markov process, then properties 1. and 2. are implied by
\begin{enumerate}
\item[3.] $X$ is locally spatial-uniformly stochastically continuous from the right: 
\begin{equation} \label{local-cont}
\forall \varepsilon>0\ \forall R>0:\ \lim_{t\downarrow 0} \sup_{|x|<R}\Prob(|X_t-x|>\varepsilon\,|\,X_0=x)=0.
\end{equation}
\end{enumerate}
\end{proposition}
\begin{proof}
Let $X$ be a $D[0,1]$ valued process. Then it has right continuous path with left limits and therefore 1. and 2. are equivalent. 

Moreover, there exists for each $\varepsilon'>0$ an $R>0$ such that $\Prob(\|X\|_\infty \geq R)<\varepsilon'.$ For $X$ being a time homogeneous Markov process and $0\leq h \leq t$ we find
\begin{equation}
\begin{split}
\Prob(|X_{t-h}-X_t|>\varepsilon) &= \int \Prob(|x -X_h|>\varepsilon\,|\,X_0 = x)\, \Prob(X_{t-h}\in dx)\\
&\leq \sup_{|x|< R} \Prob(|X_h -x|>\varepsilon\,|\,X_0 = x) + \varepsilon',
\end{split}
\end{equation}
which implies the result.
\end{proof}
In some sense Proposition \ref{stoch-continuity} suggests that it might be possible to localize condition \eqref{cond-gisko}. In fact the following is a simple consequence of Aldous result.

\begin{theorem}
\label{localized-GiSko}
Let $X^{(n)}$ be a time homogeneous strong Markov processes satisfying \eqref{global-bound} and

\begin{equation} \label{local-GiSko}
\forall \varepsilon>0\ \forall R>0:\ \lim_{h\downarrow 0}\limsup_{n\to\infty} \sup_{\substack{|x|<R\\t\leq h}}\Prob(|X^{(n)}_t-x|>\varepsilon\,|\,X_0=x)=0,
\end{equation}

then $(X^{(n)})_{n\in\N}$ is $J_1$-tight.
\end{theorem}
\begin{proof}
Let \eqref{global-bound} 
for $\varepsilon'>0$ and \eqref{local-GiSko} 
hold, and let $(\tau_n)_{n\in\N}$ be such that $\tau_n$ is a stopping time for $X^{(n)}$. Furthermore, let $\varepsilon>0$ and $(t_n)_{n\in\N}$ be a sequence in $[0,1]$ with $t_n\to 0.$
Hence 
\begin{equation}
\Prob\left(|X^{(n)}_{(\tau_n+t_n)\land 1}-X_{\tau_n}^{(n)}|>\varepsilon\right) \leq \varepsilon' + \sup_{\substack{|x|<R\\t\leq t_n}}\Prob\left(|X_{t}^{(n)}-x|>\varepsilon\,\bigg|\,X_0=x\right)
\end{equation}
and Theorem \ref{thm-aldo} implies the result.
\end{proof}

So far we have discussed conditions for Markov processes. In the following we will adapt these conditions to the Markov chain setting.

Let $Y^{(n)}$ be a time homogeneous Markov chain $(Y_0^{(n)}, Y_1^{(n)}, Y_2^{(n)},\ldots)$ on $(\Omega,\Askript,\Prob)$ and define the embeddings $X^{n,\Tau}(\omega)$ analogous to the previous section for each $\omega$. Then each $X^{n,\Tau}$ is a $D[0,1]$-valued random variable. But in general, $X^{n,\Tau}$ is not a Markov process!

Starting with a Markov chain a Markov processes can be constructed by subordination: let $(N_t)_{t\geq 0}$ be a Poisson process with intensity 1, which is independent of the Markov chain. Then one can embed the Markov chain $Y^{(n)}$ into a continuous time Markov process $\big(Z^{(n)}_t\big)_{t\geq 0}$ by setting  
\begin{equation}\label{markov-emb}
Z^{(n)}_t:=Y^{(n)}_{N_{nt}}\text{ for }t \in [0,1)\text{ and }Z^{(n)}_{1} := Z^{(n)}_{1-}.
\end{equation} The $J_1$ embedding and the Markov embedding are closely related as the following (technical) result shows.

\begin{lemma} \label{markov-vs-J1} Let 
\begin{equation}\label{extra-steps}
\sup_{k<|n-1-N_{n-}|}\left|Y^{(n)}_{((n-1)\land N_{n-} )+k} - Y^{(n)}_{(n-1)\land N_{n-} }\right| \xrightarrow{\Prob} 0.
\end{equation}
Then $X^{n,J_1}$ converges in distribution w.r.t. $J_1$ if and only if $Z^{(n)}$ converges in distribution w.r.t. $J_1.$ 
\end{lemma}
\begin{proof}
Recall that $X_t^{n,J_1} = Y^{(n)}_{\lfloor nt \rfloor}$ and $Z_t^{(n)} = Y^{(n)}_{N_{nt}}$ for $t<1.$ The first $(n-1)\land N_{n-}$ steps of these processes coincide by definition, they just appear at different times ($\frac{k}{n}$ vs. $k$-th jump time of $N_{nt}$). By a time change with a piecewise linear function $\lambda\in \Lambda_c$ both paths (up to the waiting time after the $(n-1)\land N_{n-}$-th jump) can be made to coincide. The value of $\|\lambda - id\|_\infty$ is attained at one of the jump times, thus (since $N_{n.} = \lfloor n \frac{N_{n.}}{n}\rfloor$) one can show that
\begin{equation}
\|\lambda - id\|_\infty = \sup_{s\in[0,1)}\left|s-\frac{N_{ns}}{n}\right|.
\end{equation}
The steps after $(n-1)\land N_{n-}$ can not be compensated by a time transformation. They have to be estimated explicitly. Therefore
\begin{equation}
\begin{split}
d_{J_1}\left(X_.^{n,J_1}, Z_.^{(n)}\right) \leq &\sup_{s\in[0,1)}\left|s-\frac{N_{ns}}{n}\right|\\
&+\sup_{k<|n-1-N_{n-}|}|Y^{(n)}_{((n-1)\land N_{n-} )+k} - Y^{(n)}_{((n-1)\land N_{n-} )}|.
\end{split}
\end{equation}
Since $s-\frac{N_{ns}}{n}$ is a martingale we find with Doob's maximal inequality 
\begin{equation}
\begin{split}
\E \left(\sup_{s\in[0,1)}\left|s-\frac{N_{ns}}{n}\right|\right) & \leq \sqrt{\E \left(\sup_{s\in[0,1)}\left|s-\frac{N_{ns}}{n}\right|^2\right)}\\
&\leq 2 \sqrt{\E \left(1-\frac{N_n}{n}\right)^2} =  \frac{2}{\sqrt{n}}.
\end{split}
\end{equation}
This implies with \eqref{extra-steps} that the $J_1$ distance of $X^{n,J_1}$ and $Z^{(n)}$ converges in probability to 0. Thus the convergence in distribution of either $X^{n,J_1}$ or $Z^{(n)}$ implies also the convergence in distribution of the other, e.g. by \cite[Theorem 4.1, p. 25]{Bill1968}.
\end{proof}
Before analyzing condition \eqref{extra-steps} consider the question we have asked at the beginning: when does the converse of Corollary \ref{simpleimplication} hold. Suppose a step embedded (i.e., using the $J_1$-embedding of \eqref{def-embedding}) Markov chain converges for example in $J_2$ but not in $J_1$, then the limit (before identifying it with a $D[0,1]$-function) has to have some states which it reaches by a jump and leaves instantaneously by an other jump. The following condition is sufficient to ensure that such limit points do not exist:

\begin{equation}\label{tightness} 
\forall \varepsilon>0\, \forall R>0:\ \lim_{h\to 0}\limsup_{n\to \infty}\sup_{\substack{|x|<R\\t\leq h}}\Prob\left(\left|Y_{\lfloor t n\rfloor}^{(n)}-x\right| > \varepsilon\,\bigg|\,Y_0^{(n)} = x\right)=0. 
\end{equation}
Note that this is the Markov chain version of \eqref{local-GiSko}. It ensures, as \eqref{local-GiSko}, that the limit process is locally spatial-uniformly stochastically continuous from the right and together with the Markov chain version of \eqref{global-bound}, i.e.,
\begin{equation}
\label{MC-global-bound}
\forall \varepsilon>0\,\forall m\in\N\, \exists R>0: \ \sup_n \Prob\left(\|Y_{\lfloor \cdot nm\rfloor}^{(n)}\|_\infty > R\right) <\varepsilon 
\end{equation} 
we will get a $J_1$-tightness condition, see Theorem \ref{J1-mc-convergence}. Note that in \eqref{MC-global-bound} the extra $m$ is needed since the subordinated chain might have more than $n$ steps. In the context of $D[0,1]$ this might seem surprising, but for processes in $D[0,\infty)$ the condition remains unchanged and becomes natural (cf. the last paragraph of this section). Now we can relate \eqref{extra-steps} to these conditions.

\begin{lemma}
\label{bound-extra-steps}
Let \eqref{tightness} and \eqref{MC-global-bound} hold. Then for any sequence $\left(l(n)\right)_{n\in\N}\subset [0,\infty)$ with $\frac{l(n)}{n}\to 0$
\begin{equation}
\label{exittime-bound}
\lim_{n\to \infty} \sup_{|x|<R} \Prob\left(\sup_{k\leq l(n)} |Y_k^{(n)}-x|>\varepsilon\,\bigg|\,Y_0^{(n)}=x\right) = 0
\end{equation}
and \eqref{extra-steps} hold. 
\end{lemma}

\begin{proof} Let $\tau_{B_\varepsilon(x)}^{(n)}$ denote the time of the first exit of $X^{n,J_1}$ from the ball with center $x$ and radius $\varepsilon$. Then \eqref{exittime-bound} becomes
\begin{equation}\label{exittime-explicit}
\lim_{n\to\infty} \sup_{|x|<R} \Prob\left(\tau_{B_\varepsilon(x)}^{(n)} \leq \frac{l(n)}{n} \, \bigg|\, X^{n,J_1}_0=x\right) = 0.
\end{equation}
Suppose that the limit in \eqref{exittime-explicit} is not zero. Then the limiting process (if it exists) would not be locally spatial-uniformly stochastically continuous from the right, and this contradicts \eqref{tightness}. Hence \eqref{exittime-bound} must hold. Alternatively, for a direct proof note that analogous to \cite[Lemma 2, p. 420]{GikhSkor1974} one gets for $n$ large
\begin{equation}
\begin{split}
&\Prob\left(\sup_{k\leq l(n)} |Y_k^{(n)}-x|>\varepsilon \,\bigg|\, Y_0^{(n)} =x\right)\\ 
& \leq \frac{\Prob\left(|Y_{l(n)}^{(n)}-x|>\frac{\varepsilon}{2}\,\bigg|\,Y_0^{(n)}=x\right)}{1-\sup_{\substack{|y|<R\\k\leq l(n)}}\Prob\left(|Y^{(n)}_k-y|>\frac{\varepsilon}{2}\,\bigg|\,Y_0^{(n)}=y\right)-\varepsilon_R}
\end{split}
\end{equation}
where $\varepsilon_R\in[0,1)$ is some constant depending on $R$. Hence the statement follows by \eqref{tightness} and the estimate
\begin{equation}
\begin{split}
&\limsup_{n\to\infty} \sup_{|x|<R} \Prob\left(|Y_{l(n)}^{(n)} - x| > \textstyle\frac{\varepsilon}{2}\,\bigg|\, Y_0^{(n)} =x\right) \\
& \leq \limsup_{n\to\infty} \sup_{\substack{|x|<R\\t\leq\frac{l(n)}{n}}} \Prob\left(|Y_{\lfloor tn\rfloor}^{(n)} - x| > \textstyle\frac{\varepsilon}{2}\,\bigg|\, Y_0^{(n)} =x\right)\\
& \leq \limsup_{n\to\infty} \sup_{\substack{|x|<R\\t\leq h}} \Prob\left(|Y_{\lfloor tn\rfloor}^{(n)} - x| > \textstyle\frac{\varepsilon}{2}\,\bigg|\, Y_0^{(n)} =x\right)
\end{split}
\end{equation}
which holds for any $h\in(0,1]$.

For the second part of the statement let $\varepsilon,\varepsilon' >0$, $m\in\N$ and, using \eqref{MC-global-bound}, $R$ such that $\Prob(\|Y_{\lfloor \cdot n(m+2)\rfloor}\|_\infty > R) < \varepsilon'$. Note that $\Prob(|n-1-N_{n-}|\geq \sqrt{n} m) \leq \frac{1}{m}$. Thus
\begin{equation}
\begin{split}
&\Prob\left(\sup_{k<|n-1-N_{n-}|}\left|Y^{(n)}_{((n-1)\land N_{n-} )+k} - Y^{(n)}_{(n-1)\land N_{n-} }\right|>\varepsilon\right) \\
\leq\  &\varepsilon' + {\textstyle\frac{1}{m}} + \sup_{|x|<R} \Prob\left(\sup_{k\leq \sqrt{n}m} |Y_k^{(n)}-x|>\varepsilon\,\bigg|\,Y_0^{(n)}=x\right)
\end{split}
\end{equation}
and \eqref{exittime-bound} implies \eqref{extra-steps}, since $\frac{\sqrt{n}m}{n} \to 0.$
\end{proof}

Now we can prove a $J_1$-tightness condition for embedded Markov chains, i.e., conditions 2. and 3. in the following Theorem.

\begin{theorem} \label{J1-mc-convergence}
Let $Y^{(n)},$ $X^{n,J_1}$ and $X$ be as above. Suppose the following conditions hold:
\begin{enumerate} 
\item $X^{n,J_1} \xrightarrow{\text{fdd}} X$ on a dense subset of $[0,1],$
\item \eqref{MC-global-bound}, i.e., $$\forall \varepsilon>0\,\forall m\in\N\, \exists R>0: \ \sup_n \Prob(\|Y_{\lfloor \cdot nm\rfloor}^{(n)}\|_\infty > R) <\varepsilon,$$ 
\item \eqref{tightness}, i.e., $$\forall \varepsilon>0\, \forall R>0:\ \lim_{h\to 0}\limsup_{n\to \infty}\sup_{\substack{|x|<R\\t\leq h}}\Prob\left(\left|Y_{\lfloor t n\rfloor}^{(n)}-x\right| > \varepsilon\,\bigg|\,Y_0^{(n)} = x\right)=0.$$ 
\end{enumerate}
Then
$$X^{n,J_1}\xrightarrow{d} X \text{ w.r.t. } J_1.$$ 
\end{theorem}
\begin{proof}
Assume that the conditions hold and let $Z^{(n)}$ be the Markov embedding of $Y^{(n)}$ as defined in \eqref{markov-emb}. By Lemma \ref{markov-vs-J1} and Lemma \ref{bound-extra-steps}
\begin{equation}
X^{n,J_1} \xrightarrow{d} X \text{ w.r.t. } J_1 \quad \Leftrightarrow \quad  Z^{(n)}\xrightarrow{d} X \text{ w.r.t.  } J_1.
\end{equation}
Next note that 
\begin{equation}\label{eq-poisson-trick}
\begin{split}
\Prob\left(\|Y_{N_{.n}}^{(n)}\|_\infty > R\right) &= \Prob\left(\sup_{k\leq N_n}|Y_k^{(n)}|>R\right) \\
&= \sum_{l=1}^\infty \Prob\left(\sup_{k\leq l}|Y_k^{(n)}|>R\right)\Prob(N_n = l)\\
& \leq \Prob\left(\sup_{k\leq nm} |Y_k^{(n)}|>R\right) + \Prob(N_n\geq mn)
\end{split}
\end{equation}
with $\Prob(N_n\geq mn)\leq \frac{1}{m}$ and condition 2. implies that $Z^{(n)}$ satisfies \eqref{global-bound}. Furthermore, let $\varepsilon>0.$ Then as in \eqref{eq-poisson-trick}
\begin{equation}
\begin{split}
&\Prob\left(|Y_{N_{tn}}^{(n)}-x|>\varepsilon\,|\, Y_0^{(n)} =x \right) \\&= \sum_{l=0}^\infty \Prob\left(|Y_l^{(n)}-x|>\varepsilon\,|\, Y_0^{(n)} = x\right) \Prob(N_{tn} = l)\\
& \leq \sup_{l\leq \lfloor tnm\rfloor} \Prob\left(|Y_l^{(n)}-x|>\varepsilon\, |\, Y_0^{(n)} = x\right) + \frac{1}{m}.
\end{split}
\end{equation}

Condition 3. and the arbitrary choice of $m$ imply that $Z^{(n)}$ satisfies \eqref{local-GiSko}, letting therein $h<\frac{1}{m}$. Thus by Theorem \ref{localized-GiSko} the family $(Z^{(n)})_{n\in\N}$ is $J_1$-tight.

Hence, for every sequence $n_k\nearrow \infty$ there is a subsequence $n_{k_l}$ such that $Z^{(n_{k_l})}$ converges in distribution w.r.t. $J_1$ to some limit, and $X^{(n_{k_l})}$ must have the same limit in distribution. But by 1. the limit of $X^{(n_{k_l})}$ is $X$ and it is independent of the sequence. Thus $X^{n,J_1}\xrightarrow{d} X \text{ w.r.t. } J_1.$
\end{proof}

Furthermore \eqref{tightness} also yields a statement about the convergence of finite dimensional distributions when switching the embedding.

\begin{lemma} \label{lem-fdd}
Let \eqref{tightness} hold and $X^{n,\Tau} \xrightarrow{d} X$ w.r.t. $\Tau$, hence (by Theorem \ref{thm-fdd-d}) also $X^{n,\Tau} \xrightarrow{\text{fdd}} X$ on some $I$ which is a dense subset of $[0,1].$
Then 
$$X^{n,J_1} \xrightarrow{\text{fdd}} X \text{ on } I.$$
\end{lemma}
\begin{proof}
Let the assumptions hold and note that by definition of the embeddings $X_t^{n,\Tau} \in [[Y_{\lfloor tn\rfloor}^{(n)},Y_{\lfloor tn\rfloor +1}^{(n)}]].$ Fix $\varepsilon>0$ then by Remark \ref{rem-fdd-d}.4. there exists an $R>0$ such that 
\begin{equation}
\sup_{n}\Prob\left(\|Y^{(n)}_{\lfloor\cdot n\rfloor}\|_\infty \geq R\right) <\varepsilon.
\end{equation}  
Thus
\begin{equation}
\begin{split}
&\Prob\left(|X_t^{n,J_1} - X_t^{n,\Tau}|>\varepsilon\right) \\
&= \Prob\left(|Y_{\lfloor tn\rfloor}^{(n)} - X_t^{n,\Tau}|>\varepsilon\right)\\
&\leq   \Prob\left(|Y_{\lfloor tn\rfloor}^{(n)} - Y_{\lfloor tn\rfloor+1}^{(n)}|>\varepsilon\right)\\
&\leq \sup_{x< R} \Prob\left(|Y_1^{n}-x|>\varepsilon \,|\, Y_0 = x\right) + \sup_{n}\Prob\left(\|Y^{(n)}_{\lfloor \cdot n\rfloor} \|_\infty \geq R\right)
\end{split}
\end{equation}
and the first summand converges by \eqref{tightness} to $0$ as $n\to 0.$ Since $\varepsilon$ is arbitrary the result follows by Slutsky's theorem as in Lemma \ref{markov-vs-J1}.
\end{proof}

Finally we get the following extension to Corollary \ref{simpleimplication}.

\begin{corollary} \label{conv-tight-equi} Let \eqref{tightness} and \eqref{MC-global-bound} hold. Then $X^{n,\Tau}$ converges in distribution w.r.t. $\Tau$ if and only if $X^{n,J_1}$ converges in distribution w.r.t.  $J_1$.
\end{corollary}
\begin{proof}
The direction from $J_1$ to the other topologies is just Theorem \ref{embedding-equiv} and Corollary \ref{simpleimplication}.

For the converse let \eqref{tightness}, \eqref{MC-global-bound} hold and $X^{n,\Tau} \xrightarrow{d} X$ w.r.t. $\Tau$. By Lemma \ref{lem-fdd} we get $X^{n,J_1}\xrightarrow{\text{fdd}} X$ on some dense subset of $[0,1]$ and hence (with Remark \ref{rem-fdd-d}.4.) all conditions of Theorem \ref{J1-mc-convergence} are satisfied. Therefore $X^{n,J_1}\xrightarrow{d} X$ w.r.t. $J_1.$ 
\end{proof}

\begin{example} \label{rem-feller}
If $(Y^{(n)})_{n\in\N}$ is the Markov chain approximation to a Feller processes with symbol $(x,\xi)\mapsto q(x,\xi)$ (see \cite{BoetSchi2009,BoetSchiWang2013} for the definitions and further details) then $(Y^{(n)}_{N_{nt}})_{t\geq 0}$ is a Feller process with symbol $(x,\xi)\mapsto n(e^{-\frac{1}{n} q(x,\xi)}-1)$. Thus by \cite[Corollary 5.2., p. 114]{BoetSchiWang2013} exists a constant $c$ such that 
\begin{equation}
\Prob\left(\sup_{s\leq t} \left|Y^{(n)}_{N_{ns}} - x\right| >r\, \bigg|\,Y_0^{(n)} = x\right) \leq ct \sup_{|y-x|\leq r} \sup_{|\xi|\leq \frac{1}{r}} \left| n (e^{-\frac{1}{n} q(x,\xi)}-1)\right| 
\end{equation}
for all $x\in\R^d$, $r,t>0$. Hence for any $\varepsilon >0$ 
\begin{equation}
\begin{split} \label{feller-ball-exit}
&\limsup_{n\to \infty} \Prob\left(\sup_{s\leq t} \left|Y_{\lfloor n s\rfloor}^{(n)}-x\right|> r \,\bigg|\, Y_0^{(n)}=x\right) \\
&\leq \limsup_{n\to \infty} \Prob\left(\sup_{s\leq t+\varepsilon} \left|Y_{N_{ns}}^{(n)}-x\right|> r,\ N_{n(t+\varepsilon)}\geq nt \,\bigg|\, Y_0^{(n)}=x\right) \\
&\phantom{\leq \ }+ \limsup_{n\to \infty}\Prob\left(N_{n(t+\varepsilon)} < nt\right) \\
&\leq c (t+\varepsilon) \sup_{|y-x|\leq r}\sup_{|\xi|\leq \frac{1}{r}} |q(y,\xi)|.
\end{split}
\end{equation}
Now \eqref{tightness} is satisfied if the supremum is finite, e.g. for $q$ continuous. Assuming
\begin{equation}
\lim_{r \to \infty} \sup_{|y-x|\leq r}\sup_{|\xi|\leq \frac{1}{r}} |q(y,\xi)| = 0 \ \text{ for all }x\in\R^d
\end{equation}
and using \eqref{feller-ball-exit} and \cite[(5.1)]{Boet2011a} one finds an $R$ such that condition \eqref{MC-global-bound} holds for large values of $n$, say $n\geq N$. Taking the maximum of this and the finite number of $R$'s corresponding to $n<N$ yields \eqref{MC-global-bound}.

Thus we obtained a new proof of the convergence of the Markov chain approximation of Feller processes in $J_1$. Moreover, using the introduced embeddings, the approximations converge in the four Skorokhod topologies. 
\end{example}

Finally we want to emphasize that condition \eqref{tightness} ensures the local spatial-uniform stochastic continuity from the right for the limit. Hence by Proposition \ref{stoch-continuity} the limit is stochastically continuous and therefore by Remark \ref{rem-fdd-d}.3. we can consider the convergence of the finite dimensional distributions on the whole interval $[0,1]$ (no exceptional times!). Hence, if  $X^{(n)},X$ are $D[0,\infty)$ valued processes and $X$ is stochastically continuous, then the extension from $D[0,T]$ to $D[0,\infty)$ (cf. Remark \ref{extension-of-01}.2.) does not need a restriction of the time set, i.e., we have (see \cite[Theorem 3']{Lind1973}):
\begin{equation}
\begin{split}
&X^{(n)}\xrightarrow{d}X \text{ w.r.t. } J_1 \\
\Leftrightarrow\quad  &X^{(n)}\big|_{[0,T]} \xrightarrow{d} X\big|_{[0,T]} \text{ w.r.t. $J_1$ in $D[0,T]$ for all } T\in(0,\infty).
\end{split}
\end{equation}
Additionally, in this setting $X^{(n)}\xrightarrow{d}X \text{ w.r.t. } \Tau$ implies, cf. Remark \ref{rem-fdd-d}.4., that \eqref{MC-global-bound} is satisfied.

\section{Proof  of Lemma \ref{osc-relations}} \label{sec:proofs}
We start with two elementary inequalities. Let $a,b,c,d \in \R$. Then
\begin{equation} \label{eq:minsum}
(a+b) \land (c+d) \leq (a+b\lor d) \land (c + b\lor d) = a\land c + b\lor d.
\end{equation}
For $x,y,z \in \R^d$ let $\gamma \in [[x,z]]$ such that $|y-\gamma| = |y-[[x,z]]|$, then
\begin{equation} \label{eq:tri}
|x-y| \leq |x-\gamma| + |y- \gamma| \leq |x-z| + |y-[[x,z]]|.
\end{equation}

Now let $\delta>0$ and $f:[0,1] \to \R^d$. Furthermore, let $(t,t_1,t_2) \in T_1(\delta)$ and $t_1^\star \leq t_1, t_2^\star \geq t_2$ such that $(t,t^\star_1,t^\star_2) \in T_2(\delta).$ Applying \eqref{eq:tri} and \eqref{eq:minsum} yields
\begin{equation}
\begin{split}
|f(t)-&f(t_1)|\land |f(t)-f(t_2)|\\
&
\leq \phantom{\land} \bigg(\big|f(t)-f(t_1^\star)\big| + \big|f(t_1) - [[f(t),f(t_1^\star)]]\big|\bigg) \\
& \phantom{\leq} \land \bigg(\big|f(t)-f(t_2^\star)\big| + \big|f(t_2) - [[f(t),f(t_2^\star)]]\big|\bigg)\\
&\leq \phantom{+} \bigg(\big|f(t)-f(t_1^\star)\big|\land \big|f(t)-f(t_2^\star)\big|\bigg)\\
&\phantom{\leq} + \bigg( \big|f(t_1) - [[f(t),f(t_1^\star)]]\big| \lor \big|f(t_2) - [[f(t),f(t_2^\star)]]\big|\bigg)\\
&\leq \Delta_{J_2}(\delta,f) + \Delta_{M_1}(\delta,f).
\end{split}
\end{equation}
Thus
\begin{equation}
\Delta_{J_1}(\delta,f) \leq \Delta_{J_2}(\delta,f) + \Delta_{M_1}(\delta,f),
\end{equation}
since $(t,t_1,t_2) \in T_1(\delta)$ was arbitrary. \qed

\section*{Acknowledgements}
The author is grateful for the comments of an anonymous referee who helped to improve the paper.

\end{document}